\documentclass[11pt]{article}

\usepackage{enumitem}

\setlist[enumerate]{label=(\roman*)}

\usepackage{tabularx,url} %tables, e.g. in list of symbols
\usepackage{amssymb}
\usepackage{amsmath}
\usepackage{amsthm}
\usepackage{latexsym}
\usepackage{amsfonts}
\usepackage{geometry}
\allowdisplaybreaks

\usepackage{bbm} %for writing indicator functions as doule 1 with \mathbbm{1}

\numberwithin{equation}{section} %numbering equations with a section number, requires amsmath
\allowdisplaybreaks %allow breaking equations between pages

\usepackage{xcolor} %colourful mathematics

\newtheorem{Theorem}{Theorem}[section]
\theoremstyle{definition}
\newtheorem{Lemma}[Theorem]{Lemma}
\newtheorem{Definition}[Theorem]{Definition}

\newtheorem{Example}[Theorem]{Example}

\newtheorem{Remark}[Theorem]{Remark}
\newtheorem{Proposition}[Theorem]{Proposition}
\newtheorem*{Th*}{Theorem}
\newtheorem*{Lem*}{Lemma}
\newtheorem*{Def*}{Definition}
\newtheorem*{Not*}{Notation}
\newtheorem*{Fac*}{Fact}
\newtheorem*{Exa*}{Example}
\newtheorem*{Obs*}{Observation}
\newtheorem*{Cor*}{Corollary}
\newtheorem*{Rem*}{Remark}
\newtheorem*{Pro*}{Proposition}
\newtheorem*{Exe*}{Exercise}
\newtheorem*{Que*}{Question}
\newtheorem*{Prob*}{Problem}

\usepackage{mathtools}
\DeclarePairedDelimiter\abs{\lvert}{\rvert}%
\DeclarePairedDelimiter\norm{\lVert}{\rVert}%

\makeatletter
\let\oldabs\abs
\def\abs{\@ifstar{\oldabs}{\oldabs*}}
\let\oldnorm\norm
\def\norm{\@ifstar{\oldnorm}{\oldnorm*}}
\makeatother
\DeclarePairedDelimiter\bignorm{\big\Vert}{\big\Vert}

\usepackage{mathrsfs}

\newcommand{\R}{\mathbb{R}}
\newcommand{\Rp}{\mathbb{R}_{+}}

\newcommand{\N}{\mathbb{N}}
\newcommand{\C}{\mathbb{C}}

\newcommand{\ud}{\mathrm{d}}

\newcommand{\dx}{\, \mathrm{d}x}

\newcommand{\ds}{\, \mathrm{d}s}
\newcommand{\dt}{\, \mathrm{d}t}

\newcommand{\dL}{\, \mathrm{d}L}

\newcommand{\dX}{\, \mathrm{d}X}

\newcommand{\E}{\mathbb{E}}

\newcommand{\1}{\mathbbm{1}} %requires package bbm

\newcommand{\pairing}[2]{{_{V^*}\langle #1 , #2 \rangle_V}}
\newcommand{\scalar}[2]{\langle #1 , #2 \rangle}
\newcommand{\scalarH}[2]{\langle #1 , #2 \rangle_H}

\newcommand{\scalarU}[2]{\langle #1 , #2 \rangle_U}
\newcommand{\HSnorm}[1]{\norm{#1}_{L_{\rm HS}(U,H)}}
\newcommand{\bigHSnorm}[1]{\big\Vert #1 \big\Vert_{L_{\rm HS}(U,H)}}

\newcommand{\normH}[1]{\norm{#1}_H}
\newcommand{\normV}[1]{\norm{#1}_V}
\newcommand{\normU}[1]{\norm{#1}_U}
\newcommand{\normVstar}[1]{\norm{#1}_{V^*}}
\DeclareMathOperator{\Span}{\mathrm{Span}}

\DeclareMathOperator{\sgn}{\mathrm{sgn}}
\DeclareMathOperator{\Leb}{\mathrm{Leb}}
\renewcommand{\leq}{\leqslant}
\renewcommand{\geq}{\geqslant}
\renewcommand{\le}{\leqslant}
\renewcommand{\ge}{\geqslant}

\let\phi\varphi
\let\epsilon\varepsilon

\newcommand{\itemEq}[1]{%
\begingroup%
\setlength{\abovedisplayskip}{0pt}%
\setlength{\belowdisplayskip}{0pt}%
\parbox[c]{\linewidth}{\begin{flalign}#1&&\end{flalign}}%
\endgroup}

\DeclareMathOperator{\Borel}{{\mathfrak B}}

\title{Variational solutions of stochastic partial differential equations with cylindrical L\'evy noise}
\author{Tomasz Kosmala \thanks{tomasz.kosmala@kcl.ac.uk} \qquad\qquad\qquad  Markus Riedle \thanks{markus.riedle@kcl.ac.uk} \\Department of Mathematics \\ King's College London\\ London WC2R 2LS\\ United Kingdom}
\date{16 December 2019}

\begin{document}

\maketitle

\begin{abstract}
In this article, the existence of a unique  solution in the variational approach of the stochastic evolution equation 
$$\dX(t) = F(X(t)) \dt + G(X(t)) \dL(t)$$
driven by a cylindrical L\'evy process $L$   is established. The coefficients $F$ and $G$ are assumed to satisfy the usual  monotonicity and coercivity conditions. The noise is modelled by a cylindrical L\'evy processes which is assumed to belong to a certain subclass of cylindrical L\'evy processes and may not have finite moments. 
\end{abstract}
{\bf Keywords and Phrases:} cylindrical L\'evy processes, stochastic partial differential equations, multiplicative noise, variational solutions. 
% 	60G51   	Processes with independent increments; Lévy processes
% 	60H15   	Stochastic partial differential equations 
% 	60G20   	Generalized stochastic processes
% 	28C20   	Set functions and measures and integrals in infinite-dimensional spaces (Wiener measure, Gaussian measure, etc.)
%   46G12   	Measures and integration on abstract linear spaces
%   46T12   	Measure (Gaussian, cylindrical, etc.) and integrals (Feynman, path, Fresnel, etc.) on manifolds

%\newpage
%\tableofcontents

\section{Introduction}
\label{sec_introduction}

Cylindrical Brownian motion, or equivalently Gaussian space-time noise, has been the standard model for the driving noise for stochastic partial differential equations (SPDEs) for the last 50 years. Cylindrical Brownian motions are naturally extended to the class of cylindrical L\'evy processes in order to capture non-Gaussian and discontinuous random perturbations; see Applebaum and Riedle \cite{Applebaum_Riedle}.  However, in contrast to genuine L\'evy processes, cylindrical L\'evy processes do not enjoy a semimartingale decomposition and do not attain values in the underlying space. Consequently, conventional methods, such as martingale arguments or stopping time techniques, are no longer applicable in the general case. 

If the driving noise is only additive, various specific examples of cylindrical L\'evy processes as driving noise of (linear) SPDEs are considered; see e.g.\ Brze\'zniak and Zabczyk \cite{Brzezniak_Zabczyk}, Peszat and Zabczyk  \cite{Peszat_Zabczyk_time_regularity} and Priola and Zabczyk \cite{Priola_Zabczyk_1st}. A general framework of linear equations with additive noise modelled by arbitrary cylindrical L\'evy processes is developed in Riedle \cite{Riedle_OU} and Kumar and Riedle \cite{Kumar_Riedle-2} and \cite{Kumar_Riedle}. The case of an SPDE with a multiplicative perturbation is only considered in Riedle \cite{Riedle_L2}, however under the restrictive assumption of weak square-integrability of the driving cylindrical L\'evy process. In this case, similar methods as for cylindrical Brownian motions can be exploited. 

It is the first time in this work, that an SPDE with multiplicative noise is investigated, where the random perturbation is modelled by a cylindrical L\'evy process without assuming any conditions on the moments.
The existence of a solution in this situation cannot be necessarily anticipated from the square-integrable case or other results: it is known that the irregular jumps of a cylindrical L\'evy process, in particular in the case without moments, can cause completely novel phenomena; see e.g.\ Brze\'zniak et al \cite{Brzezniak_et_al} and  Priola and Zabczyk  \cite{Priola_Zabczyk_1st}.  We focus on a subclass of cylindrical L\'evy processes which are extensively investigated in the literature as driving noise of additive equations; see e.g.\ \cite{Brzezniak_et_al, Liu_Zhai_first, Peszat_Zabczyk_time_regularity, Priola_Zabczyk_2nd, Priola_Zabczyk_1st}.

More specifically, in this paper we consider  an evolution equation of the form
\begin{equation}\label{SPDE_introduction}
\dX(t) = F(X(t)) \dt + G(X(t)) \dL(t),
\end{equation}
in the variational approach where $L$ is a cylindrical L{\'e}vy process. The coefficients $F$ and $G$ are given operators and are assumed to satisfy standard assumptions such as monotonicity and coercivity. The variational approach, first in a deterministic and then in a stochastic setting, goes back to the works by Bensoussan, Lions and Pardoux; a brief history of the approach can be found in \cite{Krylov_Rozovskii}. Existence results for equations of the form \eqref{SPDE_introduction} but driven by a Brownian motion were derived in Krylov and Rozovskii \cite{Krylov_Rozovskii}. In a series of publications \cite{Gyongy, Gyongy_Krylov_I, Gyongy_Krylov_II}, Gy\"ongy and Krylov   generalised these results to semimartingales as driving noises. The variational approach has been extended in many directions. It is especially worth to mention the works of Liu and R\"ockner \cite{Liu_Rockner_book} and  Pr\'ev\^{o}t and R\"ockner \cite{Prevot_Rockner}, where the assumptions on the coefficients were relaxed, such that classical models for example from fluid dynamics are captured by the framework. Recently, Brze\'zniak, Liu and Zhu \cite{Brzezniak_Liu_Zhu} have considered  equations  of the form \eqref{SPDE_introduction} with locally monotone coefficients and driven by a L\'evy-type noise, i.e.\ a noise which originates from a generalisation of the Gaussian space-time white noise.

The driving noise $L$ in \eqref{SPDE_introduction} is assumed to be a member of a certain subclass of 
cylindrical L\'evy processes. Elements of this subclass may not have finite moments but they enjoy a  certain kind of semimartingale decomposition. Nevertheless, classical stopping time arguments or other standard techniques cannot be exploited since 
the unbounded potential of the cylindrical  L\'evy process causes divergence of the Galerkin approximation. We circumvent this problem by introducing weights in each dimension accordingly. We present several examples of cylindrical L\'evy processes of the subclass under consideration, and relate those to models considered in the literature.

The organisation of the paper is as follows. In Section \ref{sec_preliminaries}, we collect some preliminaries and give the definition of cylindrical L\'evy processes.  
In Section \ref{sec_finite_second}, we establish the existence of a solution of equation \eqref{SPDE_introduction} in the case of a cylindrical L\'evy process with weak second moments. Since the arguments are similar to the classical case, investigated in  Brze\'zniak, Liu and Zhu \cite{Brzezniak_Liu_Zhu}, we omit details, which can be found in Kosmala \cite{Tomasz-PhD}. 
In Section \ref{sec_jumps}, we introduce the subclass of cylindrical L\'evy processes under consideration, and define their stochastic integrals. 
Section \ref{sec_existence} is devoted to establishing the existence of a unique solution of \eqref{SPDE_introduction} driven by a cylindrical L\'evy process from this subclass.  

\section{Preliminaries}
\label{sec_preliminaries}

Let $U$ and $H$ be separable Hilbert spaces with norm $\normU{\cdot}$ and 
inner product $\scalarU{\cdot}{\cdot}$ and analogously for $H$. We fix an orthonormal basis $(e_j)$ of $U$ and $(f_j)$ of $H$. The Borel $\sigma$-algebra is denoted by $\Borel(U)$. 
We use $L(U,H)$ to denote the space of bounded operators from $U$ to $H$ equipped with the operator norm. The subspace of Hilbert-Schmidt operators from $U$ to $H$ is denoted by $L_{\rm HS}(U,H)$ and it is equipped with the norm
$$\HSnorm{\phi}^2 :=  \sum_{j=1}^\infty \normH{\phi e_j}^2.$$ 

Let $(S,\mathcal{S},\mu)$ be a $\sigma$-finite measure space. We denote 
by $L^p(S;U)$ the Bochner space of all equivalence classes of measurable functions $f\colon S\to U$ which are $p$-th integrable with respect to $\mu$ equipped with the usual norm.  We use $L^0(S;U)$ to denote the space of all equivalence classes of measurable functions $f\colon S\to U$ with the topology induced by convergence in measure. The underlying measure $\mu$ and the $\sigma$-algebra $\mathcal{S}$ are always obvious from the context, e.g.\ if $S=[0,T]\times \Omega$ then $\mu= \dt\otimes P$ and ${\mathcal S}=\Borel([0,T])\otimes {\mathcal F}$, where $\dt$ is the Lebesgue measure on $[0,T]$ and $(\Omega,\mathcal{F},P)$ is a probability space. 

For a subset $\Gamma$ of $U$, sets of the form 
\[ C(u_1, ... , u_n; B) :=\{ u \in U: (\langle u, u_1 \rangle, ... , \langle u, u_n \rangle) \in B\},\]
for $u_1, ... , u_n \in \Gamma$ and $B\in \Borel(\mathbb{R}^n)$ are called {\em cylindrical sets with respect to $\Gamma$}; the set of all these cylindrical sets is denoted by $\mathcal{Z}(U,\Gamma)$ and it is an algebra. If $\Gamma$ is finite then it is a $\sigma$-algebra. A function $\lambda \colon \mathcal{Z}(U,U) \to [0,\infty]$ is called a \emph{cylindrical measure}, if for each finite subset $\Gamma \subseteq U$ the restriction of $\lambda$ on the $\sigma$-algebra $\mathcal{Z}(U, \Gamma)$ is a measure. 
A cylindrical measure is called a cylindrical probability measure if $\lambda(U) =1$.
A \emph{cylindrical random variable} $Z$ in $U$ is a linear and continuous map $Z\colon U \rightarrow L^0(\Omega; \mathbb{R})$.
Each cylindrical random variable $Z$ defines a cylindrical probability measure $\lambda$ by 
\begin{align*}
\lambda\colon  \mathcal{Z}(U,U) \to [0,1],\qquad
\lambda(C)=P\big( (Zu_1,\dots, Zu_n)\in B\big),
\end{align*}
for cylindrical sets $C=C(u_1, ... , u_n; B)$. The cylindrical probability measure $\lambda$ is called the {\em cylindrical distribution} of $Z$. 
The characteristic function of a cylindrical random variable $Z$ is defined by
\[\phi_{Z}\colon U \rightarrow \mathbb{C}, \qquad \phi_{Z}(u)=\E [\exp (iZu)],\]
and it uniquely determines the cylindrical distribution of $Z$.

Let $(\mathcal{F}_t : t \geq 0)$ be a filtration satisfying the usual conditions. A family of cylindrical random variables $L(t)\colon U \to L^0(\Omega;\R)$, $t\geq 0$, is called a \emph{cylindrical L\'evy process} if for any $n\in \N$ and $u_1, \ldots ,u_n \in U$ we have that $\big((L(t)u_1,\ldots ,L(t)u_n) : t \geq 0\big)$ is a L\'evy process in $\mathbb{R}^n$ with respect to the filtration $(\mathcal{F}_t)$.
A version of this definition appeared for the first time in Applebaum and Riedle \cite{Applebaum_Riedle} with further modifications in \cite{Riedle_L2}. Here, we include a filtration in the definition.
The characteristic function of $L(t)$ can be written in the form
\begin{equation*}
\phi_{L(1)}(u)
= \exp \left( i p(u) - \tfrac12 q(u) + \int_U \left( e^{i \langle u,x \rangle} -1 -i\langle u,x \rangle \1_{B_\R}(\langle u,x \rangle) \right) \nu(\mathrm{d}x) \right);
\end{equation*}
see \cite[Th.\ 2.7]{Applebaum_Riedle} or \cite[Th.\ 3.4]{Riedle_infinitely}.
In the above formula,  $B_\R$ is the closed unit ball in $\R$,
$p\colon U \to \R$ is a continuous function with $p(0)=0$,
$q\colon U \to \R$ is a quadratic form,
and $\nu$ is a cylindrical measure on $\mathcal{Z}(U,U)$ satisfying
\begin{equation}\label{cylindrical_Levy_measure}
\int_U \left( \scalar{u}{v}^2 \wedge 1 \right) \, \nu(\mathrm{d}v) < \infty \qquad \mathrm{for\;all\;}u \in U.
\end{equation}
A cylindrical measure satisfying \eqref{cylindrical_Levy_measure} is called a \emph{cylindrical L\'evy measure}.

We say that $L$ is \emph{weakly square-integrable} or that it has \emph{weak second moments} if $\E \big[ \abs{L(t)u}^2 \big] < \infty$ for all $t\geq 0$ and $u\in U$.
In this case, it follows from the closed graph theorem that  $L(t)\colon U \to L^2(\Omega;\R)$ is continuous for each $t\geq 0$.
Similarly $L$ is said to be \emph{weakly mean-zero} if $\E \left[ L(t)u \right] = 0$ for $t\geq 0$ and $u \in U$. For a weakly mean-zero cylindrical L\'evy process $L$, the \emph{covariance operator} $Q\colon U\to U$  is a non-negative and symmetric linear operator defined by $\langle Qu,v \rangle = \E \big[ L(1)u L(1)v \big]$ for each $u,v\in U$.

\section{Case of finite second moments}
\label{sec_finite_second}

In this section $L$ is assumed to have weak second moments. Let $Q\colon U\to U$ be its covariance operator. 
We improve the theory of stochastic integration with respect to a weakly square-integrable cylindrical  L\'evy processes introduced in \cite{Riedle_L2} by extendind the space of integrands so that it depends on the characteristics of the integrator. 
Let $\mathcal{H}=Q^{1/2}U$ be equipped with the scalar product $\langle u,v \rangle_\mathcal{H} = \scalarU{Q^{-1/2}u}{Q^{-1/2}v}$ for $u,v \in \mathcal{H}$.
Then $(Q^{1/2}e_n)$ is an orthonormal basis of $\mathcal{H}$ and $\psi \in L_{\rm HS}(\mathcal{H},H))$ if and only if $\psi Q^{1/2}$ in $L_{\rm HS}(U,H)$ and their norms coincide; see \cite[Sec.\ 2.3.2]{Prevot_Rockner}.

The approach is based on the observation that the cylindrical increments of a cylindrical L\'evy process can be Radonified by a {\em random} Hilbert-Schmidt mapping. 
More specifically, for $0\le s\le t$ let $\Phi\colon \Omega\to L_{\rm HS}(\mathcal{H},H)$ be a simple, ${\mathcal F}_s$-measurable random  variable of the form 
\begin{equation*}
\Phi(\omega)  = \sum_{j=1}^{m} \1_{A_{j}}(\omega) \phi_{j}
\end{equation*}
for deterministic operators $\phi_{j}\in L_{\rm HS}(\mathcal{H},H)$ and sets $A_{j} \in \mathcal{F}_{s}$ for $j=1,\dots, m$. The Hilbert-Schmidt property implies that for each $j\in \{1,\dots, m\}$ there exists a {\em genuine} random variable $J_{s,t}\phi_j\colon \Omega\to H$ such that $\big(L(t)-L(s)\big)(\phi_{j}^\ast h)=\scalar{J_{s,t}\phi_{j}}{h}$  for all $h\in H$. 
By linearity one can define a random variable $J_{s,t}\Phi\colon \Omega\to H$ satisfying 
\begin{equation}
\label{property_of_J_s_t}
\scalar{J_{s,t}\Phi}{h}= \sum_{j=1}^m \1_{A_j}\big(L(t)-L(s)\big)(\phi_j^\ast h)
\qquad\text{for all }h\in H.
\end{equation}
By beginning with simple stochastic processes $(\Psi(t):\, t\in [0,T])$ of the form
\begin{equation*}
\Psi(t) = \Phi_0 \1_{\{0\}}(t)+ \sum_{k=1}^{N-1} \Phi_k \mathbbm{1}_{(t_k,t_{k+1}]}(t)
\qquad\text{for  }t\in [0,T],
\end{equation*}
where $0=t_1< t_2< \cdots < t_N=T$ and each $\Phi_k$ is a simple $\mathcal{F}_{t_k}$-measurable $L_{\rm HS}(U,H)$-valued random variable, one can define the stochastic integral as
$$\int_0^T \Psi(s)\dL(s) = \sum_{k=1}^{N-1} J_{t_k,t_{k+1}} \Phi_k.$$
Using the It\^{o} isometry (see \cite{Riedle_L2})
$$\E \left[ \normH{\int_0^T \Psi(s) \dL(s)}^2 \right] = \int_0^T \E \left[ \bigHSnorm{\Psi(s) Q^{1/2}}^2 \right] \ud s$$
one can extend the definition to all stochastic processes in the space
$$\Lambda
:= \bigg\{ \Psi\colon [0,T]\times\Omega \to L_{\rm HS}(\mathcal{H},H) : \text{predictable, } \E \left[ \int_0^T \bigHSnorm{\Psi(t)Q^{1/2}}^2 \dt \right] < \infty \bigg\}.$$
The space $\Lambda$ becomes a Banach space with the norm defined by 
$$\norm{\Psi}_\Lambda = \left( \E \left[ \int_0^T \bigHSnorm{\Psi(t)Q^{1/2}}^2 \dt \right] \right)^{1/2},\qquad \Psi \in \Lambda.$$

In the following theorem we summarise the properties of this integral. If $L$ is a genuine L\'evy process, formulas for the angle bracket processes are well known; see  \cite[Cor.\ 8.17]{Peszat_Zabczyk}.
Recall that for a square-integrable, $H$-valued martingale $M$, the angle bracket process  $\langle M,M \rangle$ is defined as the unique increasing, predictable process such that $\big(\norm{M(t)}^2 - \langle M,M \rangle (t) : t\geq 0\big)$ is a martingale.
\begin{Theorem}
\label{th_properties_integral}
Suppose that $L$ is weakly square-integrable cylindrical L\'evy process and $\Psi \in \Lambda.$ 
\begin{enumerate}
\item The process
$$\left(I(t) = \int_0^t \Psi(s)\dL(s):\, t\in [0,T]\right)$$
is a square integrable martingale and has a modification with c\`adl\`ag trajectories.
\item \label{it_angle_bracket} $\displaystyle \big\langle I(\Psi),I(\Psi) \big\rangle (t) = \int_0^t \bigHSnorm{\Psi(s)Q^{1/2}}^2 \ds\;$ for all $t\in [0,T]$ $P$-a.s.
\item \label{it_stable_under_stopping} For any  stopping time $\tau$ with $P(\tau \leq T) = 1$ we have
$$\int_0^{t \wedge \tau} \Psi(s) \dL(s) = \int_0^t \Psi(s) \1_{\{s \leq \tau\}} \dL(s)\qquad\text{for all } t \in [0,T] \text{ $P$-a.s.}$$
\end{enumerate}
\end{Theorem}

\begin{proof}
Similar results are discussed for integrals with respect to classical L\'evy processes, for (i) see \cite[Th. 8.7]{Peszat_Zabczyk}, for (ii) see \cite[Cor. 8.17]{Peszat_Zabczyk} and for (iii) see \cite[Lem.\ 2.3.9]{Prevot_Rockner}. The derivations for the cylindrical case using the properties of the Radonified increments such as \eqref{property_of_J_s_t} are straightforward and can be found in \cite{Tomasz-PhD}.
\end{proof}

%\begin{Remark}
%\label{rem_rewrite_as_genuine}
%As was pointed to us by the reviewer of an earlier version of this paper, there is another approach to stochastic integration with respect to a cylindrical, weakly square-integrable cylindrical L\'evy process with a diagonal covariance operator $Q$. Suppose that $(e_n)$ is an orthonormal basis consisting of eigenvectors of $Q$. Let
%$$\tilde{L}(t) = \sum_{n=1}^\infty \frac{1}{\sqrt{2^n}} (L(t)e_n) e_n$$
%and for any operator $\psi \in L_{HS}(\mathcal{H},H)$ we define
%$$\tilde{\psi}(u) = \sum_{n=1}^\infty \sqrt{2^n} \scalarU{u}{e_n} \psi(e_n).$$
%One can show that
%$$\int_0^T \Psi(s) \dL(s) = \int_0^T \tilde{\Psi}(s) \, \ud \tilde{L}(s)$$
%and the integral on the right-hand side is a standard integral with respect to a genuine Hilbert-space valued process.
%\end{Remark}

We now discuss the existence and uniqueness of solution in the square-integrable case.
Let $(V, \norm{\cdot}_V)$ be a separable reflexive Banach space and let $\left(H, \langle \cdot, \cdot \rangle_H \right)$ and $(U, \langle \cdot , \cdot \rangle_U)$ be separable Hilbert spaces. Let $V^*$ and $H^*$ denote their duals. Assume that $V$ is densely and continuously embedded into $H$. That is we have a Gelfand triple
$$V \subseteq H = H^* \subseteq V^*.$$
Further, denote with $_{V^*}\langle \cdot , \cdot \rangle_{V}$ the duality pairing of $V^*$ and $V$.
For all $h\in H$ and $v \in V$ we have $_{V^*}\langle h,v \rangle_V = \langle h,v \rangle_H$ and without loss of generality we may assume that $\normH{v} \leq \normV{v}$ for $v\in V$ and $\normVstar{h} \leq \normH{h}$ for $h\in H$.

We consider the equation
\begin{equation}
\label{SPDE}
\dX(t) = F\big(X(t)\big) \dt + G\big(X(t)\big) \dL(t) \qquad \text{ for } t \in [0,T],
\end{equation}
with the initial condition $X(0)=X_0$ for a square-integrable, $\mathcal{F}_0$-measurable random variable $X_0$. 
The driving noise is a cylindrical L\'evy process on a separable Hilbert space $U$. 
In this section we assume that $L$ is a weakly mean-zero, weakly square-integrable, cylindrical L\'evy process, i.e.\ a cylindrical martingale with covariance operator ${Q\colon U\to U}$. 
In the remainder of the paper, we assume that there exists an orthonormal basis $(e_n)$ of $U$ consisting of eigenvectors of $Q$. Note that this does not follow from the conditions on $Q$, since it needs not be compact. 
The coefficients in equation \eqref{SPDE} are given by functions ${F\colon V\to V^*}$ and ${G\colon V\to L_{\rm HS}(\mathcal{H},H)}$.
More specifically, we assume  the following in this section: there are constants $\alpha, \lambda, \beta, c>0$ such that:
\begin{enumerate}[label=(A\arabic*),series=assumptions]
\item \label{item_coercivity}(Coercivity) for all $v\in V$ we have
\begin{align*}
2\pairing{F(v)}v + \bigHSnorm{G(v)Q^{1/2}}^2 + \alpha \normV{v}^2 \leq \lambda \normH{v}^2+\beta;
\end{align*}
\item \label{item_monotonicity} (Monotonicity) for all $v_1,v_2 \in V,$ we have
\begin{align*}
2\pairing{F(v_1)-F(v_2)}{v_1-v_2} + \bigHSnorm{(G(v_1)-G(v_2))Q^{1/2}}^2 \leq \lambda \normH{v_1-v_2};
\end{align*}
\item \label{item_linear_growth} (Linear growth) $\normVstar{F(v)} \leq c(1+\normV{v} )$ for all $v \in V$;
\item \label{item_hemicontinuity} (Hemicontinuity) the mapping $\R \ni s \mapsto {\pairing{F(v_1+sv_2)}{v_3}}$ is continuous for all $v_1,v_2,v_3 \in V$.
\item \label{item_Q_diagonal} The cylindrical L\'evy process $L$ is weakly mean-zero 
and is weakly square-integrable. Its covariance operator $Q$ has eigenvectors $(e_j)$, which form an orthonormal basis of $U$. 
\end{enumerate}

Conditions of this form appear in most of the papers mentioned in the introduction.
For instance in \cite{Brzezniak_Liu_Zhu} these conditions are formulated for a sum of cylindrical Brownian motion with covariance equal to the identity (in particular, satisfying \ref{item_Q_diagonal}) and a Poisson random measure.
Later, we will consider the case of a non-integrable noise. 
%By truncating the jumps, we will reduce the problem of existence and uniqueness of solutions to the case of an equation driven by a process satisfying \ref{item_Q_diagonal}.
We will obtain the existence and uniqueness of solutions by truncating the jumps and reducing the problem to the case of an equation driven by a process satisfying \ref{item_Q_diagonal}.
We now give the definition of a solution to \eqref{SPDE}, similarly as in Pr\'ev\^{o}t and R\"ockner \cite[Def.\ 4.2.1]{Prevot_Rockner} or Brze\'zniak, Liu and Zhu \cite[Def.\ 1.1]{Brzezniak_Liu_Zhu}. Since we later consider the case of a driving noise without finite moments and thus the solution cannot be expected to have finite moments,  we do not require finite expectation of the solution.

\begin{Definition}
\label{def_variational_solution_non_L2}
A \emph{variational solution} of \eqref{SPDE} is a pair $(X,\bar{X})$ of an 
$H$-valued, c\`adl\`ag adapted process $X$ and a $V$-valued, predicable 
process $\bar{X}$ such that 
\begin{enumerate}
\item $X$ equals $\bar{X}$ $\mathrm{d}t \otimes P$-almost everywhere; 
\item  $P$-a.s. $\displaystyle \int_0^T \normV{\bar{X}(t)} \dt < \infty$;
\item \itemEq{X(t) = X_0 + \int_0^t F(\bar{X}(s))\ds + \int_0^t G(\bar{X}(s))\dL(s) \text{ for all } t\in [0,T]\,\, P\text{-a.s.}\label{integral_eq_in_def_of_variational_sol}}
\end{enumerate}
We say that the solution is \emph{pathwise unique} if any two variational solutions $(X,\bar{X})$ and $(Y,\bar{Y})$ satisfy
$$P\big(X(t)=Y(t) \text{ for all } t \in [0,T]\big) = 1.$$
\end{Definition}

\begin{Theorem}
\label{th_main_existence_result}
Under Assumptions {\upshape(A1)}-{\upshape(A5)}, equation \eqref{SPDE} has a unique variational solution $(X,\bar{X})$. Moreover, the solution satisfies
\begin{align*}
\int_0^T \E\left[\norm{\bar{X}(s)}_V^2\right]\ds<\infty. 
\end{align*}
\end{Theorem}

The proof can be obtained by repeating the arguments from \cite{Prevot_Rockner} and using the properties of the integral outlined in Theorem \ref{th_properties_integral}. The details are presented in \cite{Tomasz-PhD}. 
%Alternatively, one can use Remark \ref{rem_rewrite_as_genuine} to rewrite the coefficients $F$ and $G$ and transforming $L$ into a genuine process.
%Note that conditions \ref{item_coercivity} and \ref{item_monotonicity} do not change when formulated for $\tilde{G}$ and the covariance $\tilde{Q}$ of $\tilde{L}$.

\section{Orthogonal cylindrical L\'evy processes}
\label{sec_jumps}

In this section we consider the case of a driving noise without finite moments.  Contrary  to the classical case of a genuine L\'evy process, one cannot directly apply stopping time arguments such as in \cite[Sec.\ 9.7]{Peszat_Zabczyk} or interlacing techniques such as in  \cite[Th.\ IV.9.1]{Ikeda_Watanabe}, since the cylindrical L\'evy process does not attain values in the underlying space. 

For a bounded sequence of positive real numbers $c=(c_j) \in \ell^\infty(\R_+)$ we define the sequence of stopping times by
\begin{equation*}
%\label{definition_of_tau_n_k}
\tau_n^c(k) 
:= \inf \bigg\{ t \geq 0 : \sum_{j=1}^n \left( \Delta L(t)e_j \right)^2 c_j^2 > k^2 \bigg\}\qquad\text{for each }k>0, n\in\N.
\end{equation*}
The stopping time $\tau_n^c(k)$ can be seen as the first time, the $n$-dimensional L\'evy process $\big((L(t)(c_1e_1),\dots, L(t)(c_ne_n)):\,t\ge 0\big)$ has a jump of size larger than $k$. Since $\tau_n^c(k)$ is non-increasing in $n$, we can define another sequence of stopping times by
\begin{equation}
\label{definition_of_tau_k}
\tau^c(k) := \lim_{n \to \infty} \tau_n^c(k) \qquad \text{for } k>  0.
\end{equation}
Contrary to the the case of a genuine Hilbert space-valued L\'evy process, if the noise is cylindrical the stopping times $\tau_n^c(k)$  may accumulate at zero, 
i.e.\ $\tau^c(k)=0$ $P$-almost surely, see Remark \ref{rem_constant_impossible} below.
It will turn out that the distribution of the stopping time $\tau^c(k)$ depends on the parameter
\begin{equation}
\label{definition_of_m_c_k}
m^c(k) := \sup_{n\in \N} \nu \bigg(  \bigg\{ u\in U : \sum_{j=1}^n \langle u,e_j \rangle^2 c_j^2 > k^2 \bigg\} \bigg) \qquad \text{for } k>0,
\end{equation}
where $\nu$ is the cylindrical L\'evy measure of $L$.
If $L$ is a genuine L\'evy process in $U$ then its L\'evy measure $\nu$ is finite outside each ball around $0$ and $m^c(k)\to 0$ as $k\to\infty$. 
In the cylindrical case, the situation turns out to be rather different as Proposition \ref{pro_accumulation_and_m_c_k} shows:

\begin{Proposition}  \label{pro_accumulation_and_m_c_k}
Let $L$ be a cylindrical L\'evy process with $\tau^c(k)$ defined in \eqref{definition_of_tau_k} and $m^c$ defined in \eqref{definition_of_m_c_k} for a fixed $c \in \ell^\infty(\R_+)$.
\begin{enumerate}
\item[{\rm (1)}] We have the following dichotomy for each $k>0$:
\begin{enumerate}
\item $m^c(k)=0$ $\Leftrightarrow$ $\tau^c(k)=\infty$ $P$-a.s;
\item $m^c(k) \in (0,\infty)$ $\Leftrightarrow$ $\tau^c(k)$ is exponentially distributed with parameter $m^c(k)$;
\item $m^c(k)=\infty$ $\Leftrightarrow$ $\tau^c(k) = 0$ $P$-a.s.
\end{enumerate}
\item[{\rm (2)}] We have: $\lim\limits_{k\to\infty} m^c(k)=0$ $\Leftrightarrow$ $\displaystyle \lim_{k\to\infty}\tau^c(k)=\infty$ $P$-a.s.
\end{enumerate}
\end{Proposition}

\begin{proof} 
(1) Define the mapping
\begin{align*}
\pi_n^c:U \to U, \qquad \pi_n^c(u) = \sum_{j=1}^n c_j \scalar{u}{e_j} e_j.
\end{align*}
Note that $\tau_n^c(k)$ is the time of the first jump of size larger than $k$ of the finite dimensional L\'evy process $L_n^c$ defined by
$$L_n^c(t) = \sum_{j=1}^n c_j L(t)(e_j) e_j, \qquad t \geq 0.$$
As the L\'evy measure $\nu_n^c$ of $L_n^c$ is given by  $\nu_n^c:=\nu \circ (\pi_n^c)^{-1}$, the stopping time $\tau_n^c(k)$ is exponentially distributed with parameter
$$\lambda_n^k:= \nu_n^c\big(\{u\in U: \normU{u}>k\}\big) 
= \nu \bigg( \bigg\{ u \in U : \sum_{j=1}^n c_j^2 \scalar{u}{e_j}^2 > k^2 \bigg\} \bigg).$$

(i): the very definition implies that $m^c(k)=0$ if and only if $\lambda_n^k=0$ for all $n\in\N$. The latter is equivalent to $\tau_n^c(k)=\infty$ for all $n\in\N$. 

(ii), (iii): the characteristic function $\phi_{\tau_n^c(k)}$ of $\tau_n^c(k) $ is given by 
\begin{align*}
\phi_{\tau_n^c(k)}\colon \R\to\C,\qquad
\phi_{\tau_n^c(k)}(x)=\frac{\lambda_n^k}{\lambda_n^k-ix}.
\end{align*}
As $\lambda_n^k$ monotonically increases to $m^c(k)$ as $n\to\infty$, the characteristic function $\phi_{\tau_n^k}$ converges  to the characteristic function either of the exponential distribution with parameter $m^c(k)$ or of the Dirac measure in $0$.

For establishing (2), note that monotonicity of $k\mapsto \tau^c(k)$ yields
\begin{align*}
P\left(\lim_{k\to\infty} \tau^c(k) =\infty\right)
= P\bigg( \bigcap_{t\in \N} \bigcup_{n\in \N} \bigcap_{k\geq n} \{\tau^c(k) >t\}\bigg)
=\lim_{t\to\infty} \lim_{n\to\infty} P\Big(\tau^c(n)>t\Big).
\end{align*}
Since $P(\tau^c(k) > t) = \exp(-tm^c(k))$, this completes the proof of (2).
\end{proof}

We now focus on a special class of cylindrical L\'evy processes, which, similarly as in the case of a cylindrical Brownian motion, can be represented by a sum. 
That is we assume a form of the noise as in the Karhunen-Lo\`eve theorem with independent components but without requiring finite second moments.
%However, the sum does not converge in the underlying Hilbert space.
Let $L$ be a cylindrical L\'evy process with cylindrical L\'evy measure $\nu$ and let $(e_j)$ be an orthonormal basis of $U$.
%In the reminder of this section we consider the special class of 
L is called \emph{orthogonal} cylindrical L\'evy processes if it is of the form
\begin{align} \label{series_cylindrical_Levy}
L(t)u = \sum_{j=1}^\infty \ell_j(t) \langle u,e_j \rangle
\qquad P\text{-a.s.\ for all }u\in U,\, t\ge 0,
\end{align}
where $(\ell_j)$ is a sequence of independent, not necessarily identically distributed, one-dimensional L\'evy processes.
Denote the characteristics (with respect to the standard truncation 
function $\1_{B_\R}$) of $\ell_j$ by $(b_j,s_j,\rho_j)$ for each $j\in\N$. Lemma 4.2 in \cite{Riedle_OU} guarantees that the sum in \eqref{series_cylindrical_Levy} converges and defines a cylindrical L\'evy process if and only if the characteristic functions of $\ell_j$ are equicontinuous at $0$ and the following three conditions are satisfied for every $(\alpha_j) \in \ell^2(\R)$:
\begin{enumerate}
\item \itemEq{\sum\limits_{j=1}^\infty \1_{B_\R}(\alpha_j)\abs{\alpha_j} \abs{b_j + \int_{1<\abs{x}\leq \abs{\alpha_j}^{-1}} x \, \rho_j(\mathrm{d}x)} < \infty, \label{lemma_4_2_first_condition}}
\item \itemEq{ (s_j) \in \ell^\infty(\R), \label{lemma_4_2_second_condition}}
\item  \itemEq{\sum\limits_{j=1}^\infty {\displaystyle \int_\R} \left( \abs{\alpha_j x}^2 \wedge 1 \right) \rho_j(\mathrm{d}x) < \infty.\label{lemma_4_2_third_condition}}
\end{enumerate}
Independence of the L\'evy processes ($\ell_j)$ implies that the cylindrical L\'evy measure  of $L$ has only support in $\cup\, \Span \{e_j\}$. Consequently, the function $m^c$, defined in \eqref{definition_of_m_c_k}, reduces to
\begin{align}
\label{eq.m^c-series}
m^c(k) = \sum_{j=1}^\infty \rho_j\left( \left\{x\in\R: \abs{x}>\tfrac{k}{c_j}\right\} \right) \qquad\text{for all }k>0.
\end{align}

In general, due to the non-linearity of the truncation function, cylindrical L\'evy processes do not enjoy a type of L\'evy-It\^{o} decomposition. However, the specific construction of cylindrical L\'evy processes of the form \eqref{series_cylindrical_Levy} suggests to derive a L\'evy-It\^{o} decomposition from an appropriate decomposition of the real-valued processes $\ell_j$. More precisely, for a given sequence 
$c=(c_j)\in \ell^\infty (\Rp)$ and $k>0$ we obtain  $\ell_j(t)= p_j^{c,k}(t) + m_j^{c,k}(t) + r_j^{c,k}(t)$ for all $t\ge 0$ where 
\begin{align}
p_j^{c,k}(t) &:= \left( b_j + \int_{1<\abs{x}\leq k/c_j} x \, \rho_j(\mathrm{d}x)\right)t , \label{eq.def-p} \\
m_j^{c,k}(t) &:= \sqrt{s_j} W_j(t) + \int_{\abs{x}\leq k/c_j} x \, \tilde{N}_j(t,\mathrm{d}x), \\
r_j^{c,k}(t) &:= \int_{\abs{x}>k/c_j} x \, N_j(t,\mathrm{d}x).
\label{eq.def-q}
\end{align}
Here, the process $W_j$ is a real-valued standard Brownian motion and $N_j$ is a Poisson random measure on $[0,\infty)\times \R$ with intensity measure $\mathrm{d}t\otimes \rho_j$.

In the following lemma, we summarise the conditions on the cylindrical L\'evy process such that the stopping times $\tau^c(k)$, defined in \eqref{definition_of_tau_k}, do not accumulate at zero and such that the decomposition of $\ell_j$ leads to a decomposition of the cylindrical L\'evy process:
\begin{enumerate}[assumptions]
\item \label{item_existence_of_c} there exists a sequence $c=(c_j) \in \ell^\infty(\R_+)$ such that
\begin{enumerate}
\item \itemEq{\left( p_j^{c,k}(1)\right)_{j\in\N}\in \ell^2(\R) \text{ for each } k>0; \label{assumption_drift}}
\item \itemEq{\displaystyle \sup_{j\in\N}\int_{\abs{x}\leq k/c_j} x^2
 \, \rho_j(\mathrm{d}x)<\infty \text{ for each } k>0; \label{assumption_square_integrable}} 
\item \itemEq{\displaystyle \lim_{k\to\infty} m^c(k) = 0.\label{assumption_no_accummulation}} 
\end{enumerate}
\end{enumerate}

\begin{Remark} \label{rem_square_summable_c}
Assume that $L$ is of the form \eqref{series_cylindrical_Levy}, i.e.\ Conditions \eqref{lemma_4_2_first_condition}--\eqref{lemma_4_2_third_condition} are satisfied.  For a square summable sequence $(c_j)$, condition \eqref{lemma_4_2_third_condition} implies \eqref{assumption_no_accummulation} by \eqref{eq.m^c-series} and Lebesgue's theorem.  On the other hand, if $c_j$ is constantly equal to $1$, then Condition \eqref{lemma_4_2_third_condition} implies \eqref{assumption_square_integrable}. 
Indeed, suppose for contradiction that the sequence 
$\big( \int_{\abs{x}\leq k} x^2 \, \rho_j(\mathrm{d}x): j \in \N\big)$
is unbounded. Then there exists a sequence $(\alpha_j) \in \ell^2(\R)$ such that
$$\sum_{j=1}^\infty \alpha_j^2 \int_{\abs{x}\leq k} x^2 \, \rho_j(\mathrm{d}x) = \infty,$$
which contradicts \eqref{lemma_4_2_third_condition}.

In summary, for the assumption \ref{item_existence_of_c} to hold there must be some balance between the rate of decay of the L\'evy measures $(\rho_j)$ and the rate of convergence of the sequence $c$.
\end{Remark}

\begin{Lemma}
\label{lem_cyl_decomp}
Assume that $L$ is a cylindrical L\'evy process of the form \eqref{series_cylindrical_Levy} satisfying \ref{item_existence_of_c} 
for a sequence $c\in \ell^\infty(\Rp)$. Then $L$ can be decomposed into $L(t)=P^c_k(t)+M^c_k(t)+R^c_k(t)$ for each $t\ge 0$ and $k>0$, where $P^c_k$, $M^c_k$ and $R^c_k$ are cylindrical L\'evy processes defined by
\begin{equation*} %\label{definition_of_M_and_P}
P^c_k(t)u := \sum_{j=1}^\infty p_j^{c,k}(t) \langle u,e_j \rangle, \enskip
M^c_k(t)u := \sum_{j=1}^\infty m_j^{c,k}(t) \langle u,e_j \rangle, \enskip
R^c_k(t)u := \sum_{j=1}^\infty r_j^{c,k}(t) \langle u,e_j \rangle.
\end{equation*}
The process $M^c_k$ is a weakly square-integrable cylindrical L\'evy martingale with a diagonal covariance operator $Q_k$ 
and the stopping times $\tau^c$, defined in \eqref{definition_of_tau_k}, 
satisfy $\tau^c(k)\to\infty$ $P$-a.s.\ as $k\to\infty$. 
\end{Lemma}

\begin{proof}
We write $M^c_k(t)=X(t)+Y^c_k(t)$ for each $k>0$ with
$$X(t)u := \sum_{j=1}^\infty \sqrt{s_j} W_j(t) \langle u,e_j \rangle, \qquad 
Y^c_k(t)u := \sum_{j=1}^\infty \int_{\abs{x}\leq k/c_j} x \, \tilde{N}_j(t,\mathrm{d}x) \langle u,e_j \rangle, $$
for all $u\in U$. Since condition \eqref{lemma_4_2_second_condition} implies
\begin{align*}
\E\left[ \abs{X(t)u}^2\right]
=\sum_{j=1}^\infty \abs{s_j}\scalar{u}{e_j}^2 \le \norm{s}_\infty \norm{u}^2,
\end{align*}
we obtain that $X(t)\colon U\to L^0(\Omega;\R)$ is well defined, continuous and weakly square-integrable. 
We have
$$\E\left[\abs{Y^c_k(t)u}^2\right]
=t\sum_{j=1}^\infty\langle u,e_j \rangle^2  \int_{\abs{x}\leq k/c_j} x^2 \, \rho_j(\mathrm{d}x)
< \infty$$
by \eqref{assumption_square_integrable}.
Consequently, $Y^c_k(t)$ and thus $M^c_k(t)$ are well defined, continuous and weakly square-integrable. By \eqref{assumption_drift}, the (deterministic) process $P^c_k$ is well defined.
Since $R^c_k=L-M^c_k-P^c_k$ it follows that the series in the definition of $R^c_k$ converges and that $R^c_k(t)\colon U\to L^0(\Omega;\R)$ is continuous for all $t\geq 0$.
\end{proof}

\begin{Example}\label{ex.symmetric-stable}
\label{exa_symmetric_stable}[Two-sided stable process]
An often considered example of a process given in \eqref{series_cylindrical_Levy} is for $\ell_j = \sigma_j h_j$, 
where $h_j$ are identically distributed, symmetric $\alpha$-stable L\'evy processes 
and $\sigma_j\in\R$; see \cite{Priola_Zabczyk_2nd, Priola_Zabczyk_1st}. 
In this case,  $\ell_j$ has L\'evy measure $\rho_j = \rho \circ m_{\sigma_j}^{-1}$, where $m_{\sigma_j}\colon \R\to\R$ is given by $m_{\sigma_j}(x) = \sigma_j x$ and  $\rho(\mathrm{d}x) = \frac12 \abs{x}^{-1-\alpha} \dx$.
By \cite[Ex.\ 4.5]{Riedle_OU}, formula \eqref{series_cylindrical_Levy} defines a cylindrical L\'evy process if and only if $\sigma=(\sigma_j) \in \ell^{\frac{2\alpha}{2-\alpha}}(\R)$. Moreover, $L$ is induced by a classical process if and only if $\sigma \in \ell^\alpha(\R)$. 

We show that Assumption \ref{item_existence_of_c} is satisfied
for the sequence $(c_j)\in\ell^2(\Rp)$ defined by $c_j = \abs{\sigma_j}^{\frac{\alpha}{2-\alpha}}$. Condition \eqref{assumption_drift} is trivially satisfied because each $h_j$ has no drift and the L\'evy measure is symmetric. Since
\begin{align*}
\int_{\abs{x}\leq \frac{k}{c_j}} x^2 \, \rho_j(\mathrm{d}x)
= \sigma_j^2 \int_{\abs{x}\leq \frac{k}{\abs{c_j\sigma_j}}} x^2 \, \rho(\mathrm{d}x)
%=\sigma_j^2 \frac{2}{2-\alpha} \left[ x^{2-\alpha} \right]_0^\frac{k}{\abs{c_j\sigma_j}}
= \sigma_j^2 \frac{k^{2-\alpha}}{2-\alpha} \abs{c_j \sigma_j}^{\alpha-2}
= \frac{k^{2-\alpha}}{2-\alpha}, 
\end{align*}
Condition \eqref{assumption_square_integrable} is satisfied. 
Since $(c_j)\in\ell^2(\Rp)$ by its very definition,  
Remark \ref{rem_square_summable_c} establishes
Condition \eqref{assumption_no_accummulation}.
\end{Example}

\begin{Example}\label{ex.one-sided-stable}[One-sided stable process]
We choose $\ell_j=\sigma_j h_j$ in \eqref{series_cylindrical_Levy} with $\sigma_j\in\R$ and  $h_j$ arbitrary, strictly $\alpha$-stable L\'evy process with $\alpha\in (0,2)$ and with no negative jumps. Note, that $\alpha\neq 1$, since a 1-stable L\'evy process is strictly stable if and only if its L\'evy measure is symmetric.  The characteristic function of $h_j(1)$ is given by 
\begin{align*}
\phi_{h_j(1)}(x) =
\exp\left( -\abs{x}^\alpha \left(1-i \tan\tfrac{\pi\alpha}{2} \sgn x\right) \right),
\end{align*}
see \cite[Th.\ 14.15, Def.\ 14.16]{Sato}. 
It follows that the L\'evy process $\sigma_jh_j$ has characteristics $(b_j, 0,\rho_j)$ given by 
\begin{align*}
b_j= \sigma_j  \frac{1}{c_\alpha(1-\alpha)} \sigma_j \abs{\sigma_j}^{\alpha-1},
\qquad
\rho_j(\mathrm{d}x)=\big(\rho\circ m_{\sigma_j}^{-1}\big)(\mathrm{d}x), 
\end{align*}
where $c_\alpha = - \cos\left( \frac{\alpha \pi}{2} \right) \Gamma(\alpha)$, the function $m_{\sigma_j}\colon \R\to\R$ is defined by $m_{\sigma_j}(x) = \sigma_j x$
and $\rho(\mathrm{d}x) = \1_{(0,\infty)}(x)\frac{1}{c_\alpha}x^{-1-\alpha}\dx$.

We claim that $L$ is a cylindrical L\'evy process if and only if $\sigma \in \ell^\frac{2\alpha}{2-\alpha}(\R)$. Indeed, Condition \eqref{lemma_4_2_first_condition} reduces to
\begin{align*}
\sum_{j=1}^\infty \abs{\alpha_j} \abs{b_j + \int_{1<\abs{x}\le 1/\abs{\alpha_j}} x \, \big(\rho\circ m_{\sigma_j}^{-1}\big)(\mathrm{d}x)} 
%&= c \sum_{j=1}^\infty \abs{\alpha_j} \abs{\frac{\sigma_j}{1-\alpha} \abs{\sigma_j ^{\alpha-1} + \sigma_j  \int_\frac1{\abs{\sigma_j}}^\frac1{\abs{\alpha_j\sigma_j}} x x^{-1-\alpha} \dx} \\
&= \frac{c}{c_\alpha\abs{1-\alpha}} \sum_{j=1}^\infty \abs{\alpha_j\sigma_j}^\alpha, 
\intertext{whereas Condition \eqref{lemma_4_2_third_condition} reads as}
\sum\limits_{j=1}^\infty {\displaystyle \int_\R} \left( \abs{\alpha_j x}^2 \wedge 1 \right) \rho_j(\mathrm{d}x)
&= \frac{2c}{c_\alpha(2-\alpha)\alpha}\sum_{j=1}^\infty \abs{\alpha_j\sigma_j}^\alpha. 
\end{align*}

Assumption \ref{item_existence_of_c} is satisfied with $c_j = \abs{\sigma_j}^{\frac{\alpha}{2-\alpha}}$, since Condition \eqref{assumption_drift} can be calculated as 
\begin{align*}
\sum_{j=1}^\infty \bigg(b_j + \int_{1<\abs{x}\leq \frac{k}{c_j}} x \, \rho_j(\mathrm{d}x) \bigg)^2 
%&= \sum_{j=1}^\infty \bigg(\frac{c}{c_\alpha} \frac{1}{1-\alpha} \sigma_j\abs{\sigma_j}^{\alpha-1}+ \frac{c}{c_\alpha} \sigma_j \int_\frac{1}{\abs{\sigma_j}}^\frac{k}{\abs{c_j\sigma_j}} x \, \rho (\mathrm{d}x) \bigg)^2 \\
%&= \frac{c^2}{c_\alpha^2} \sum_{j=1}^\infty \bigg(\frac{1}{1-\alpha} \sigma_j\abs{\sigma_j}^{\alpha-1} + \sigma_j \int_\frac{1}{\abs{\sigma_j}}^\frac{k}{\abs{c_j\sigma_j}} x^{-\alpha} \dx \bigg)^2 \\
%&= \frac{c^2}{c_\alpha^2} \sum_{j=1}^\infty \left(\frac{1}{1-\alpha} \sigma_j\abs{\sigma_j}^{\alpha-1} + \sigma_j \left[ \frac{1}{1-\alpha}x^{1-\alpha} \right]_\frac{1}{\abs{\sigma_j}}^\frac{k}{\abs{c_j\sigma_j}} \right)^2 \\
%&= \left(\frac{c}{c_\alpha(1-\alpha)}\right)^2 \sum_{j=1}^\infty \left(\sigma_j\abs{\sigma_j}^{\alpha-1}+ \sigma_j\left( k^{1-\alpha} \abs{c_j\sigma_j}^{\alpha-1} - \abs{\sigma_j}^{\alpha-1} \right) \right)^2 \\ 
%&= \left(\frac{ck^{1-\alpha}}{c_\alpha(1-\alpha)} \right)^2 \sum_{j=1}^\infty c_j^{2\alpha} \abs{\sigma_j}^{2\alpha} \\
= \left(\frac{ck^{1-\alpha}}{c_\alpha(1-\alpha)} \right)^2 \sum_{j=1}^\infty \abs{\sigma_j}^\frac{2\alpha}{2-\alpha}.
\end{align*}
Conditions  \eqref{assumption_square_integrable} and \eqref{assumption_no_accummulation} follow by the same arguments as in Example \ref{exa_symmetric_stable}. 
\end{Example}

\begin{Remark}
\label{rem_constant_impossible}
In both Examples \ref{ex.symmetric-stable} and \ref{ex.one-sided-stable}, Condition 
\eqref{assumption_no_accummulation} would not be satisfied for a constant level of truncation of jumps i.e.\ with $c_j=1$ for all $j\in \N$. By introducing the weights $(c_j)$ we compensate the fact that the cylindrical distribution of $L$ is not tight, i.e.\ its mass of the span of the higher nodes decays too slowly. 
\end{Remark}

\begin{Example}\label{ex.one-sided-varying}[One-sided regularly varying tails]
Recall that a measure $\mu$ concentrated on $(0,\infty)$ is said to have regularly varying tails with index $\alpha$ if
$$\lim_{x \to \infty} \frac{\mu(\lambda x,\infty)}{\mu(x,\infty)} = \lambda^{-\alpha} \qquad \text{for all } \lambda>0;$$
see \cite{Bingham, Feller_vol_2}. We choose $\ell_j = \sigma_j h_j$ in \eqref{series_cylindrical_Levy} with a sequence of independent and identically distributed L\'evy processes $h_j$ of regularly varying tails of index $\alpha \in (0,1)\cup (1,2)$. For simplifying the calculations, we assume that the characteristic function of $h_j(1)$ is given by
\begin{align*}
\phi_{h_j(1)}(x) =
\begin{cases}
 \exp\left( \int_0^\infty \left( e^{ixy}-1 -ixy\1_{B_\R}(y)\right) \rho(\mathrm{d}y)+ixb \right), &\text{ if }\alpha \in (0,1),\\
  \exp\left( \int_0^\infty \left( e^{ixy}-1-ixy \right) \rho(\mathrm{d}y)\right), &\text{ if }\alpha \in (1,2),
\end{cases}
\end{align*}
for a constant $b\in\R$. The L\'evy measure $\rho$ of $h_j$ has regularly varying tails according to \cite{Embrechts_Goldie}.

We show that if $(\sigma_j) \in \ell^{\frac{2\delta}{2-\delta}}(\R)$ for some $\delta<\alpha$, then \eqref{series_cylindrical_Levy} defines a cylindrical L\'evy process. For this purpose, we define 
\begin{equation*}
V_\delta (x) := \int_x^\infty y^\delta \1_{(1,\infty)}(y)\, \rho(\mathrm{d}y), \qquad 
U_2(x) := \int_1^x y^2 \, \rho(\mathrm{d}y), \qquad \text{for } x>0.
\end{equation*}
It follows from \cite[Prop.\ 4.2.1]{Samorodnitsky} that $V_\delta(0)<\infty$ and $U_2(\infty)=\infty$. Theorem~VII.9.2 in \cite{Feller_vol_2} implies
$$\lim_{x\to \infty} \frac{x^{2-\delta} V_\delta(x)}{U_2(x)} = \frac{2-\alpha}{\alpha-\delta}=:c,$$
and therefore there exists $M>0$ such that 
$$U_2(x) \leq \frac{2x^{2-\delta}V_\delta(x)}{c}\qquad\text{for all } x\ge M.$$
Since both $(\alpha_j)$ and $(\sigma_j)$ tend to $0$ we can assume without loss of generality that $\frac1{\abs{\alpha_j\sigma_j}}>M$ for all $j\in \N$.
For verifying Condition \eqref{lemma_4_2_third_condition} we obtain 
\begin{align*}
&\sum_{j=1}^\infty \left( \alpha_j^2 \sigma_j^2 \int_{0\le x\leq 1} x^2 \, \rho(\mathrm{d}x) + \alpha_j^2 \sigma_j^2 \int_{1<x\leq \frac1{\abs{\alpha_j \sigma_j}}} x^2 \, \rho(\mathrm{d}x) \right) \\
&\qquad = \sum_{j=1}^\infty \alpha_j^2 \sigma_j^2 \int_{0\le x\leq 1} x^2 \, \rho(\mathrm{d}x) + \sum_{j=1}^\infty \alpha_j^2 \sigma_j^2 U_2\left(\frac1{\abs{\alpha_j\sigma_j}}\right)\\
&\qquad \le 
\sum_{j=1}^\infty \alpha_j^2 \sigma_j^2 \int_{0\le x\leq 1} x^2 \, \rho(\mathrm{d}x) + \frac2{c} \sum_{j=1}^\infty \abs{\alpha_j \sigma_j}^{\delta } V_\delta\left(\frac1{\abs{\alpha_j \sigma_j}}\right).
\end{align*}
Both sums are finite because of the summability assumptions on $\alpha$ and $\sigma$.
Similarly, we derive that 
\begin{align*}
\sum_{j=1}^\infty \rho\left(  \Big[\tfrac{1}{\abs{\alpha_j\sigma_j}}, \infty \Big)\right)<\infty,
\end{align*}
which shows Condition \eqref{lemma_4_2_third_condition}. 

Similarly to the stable case in Example \ref{ex.one-sided-stable}, the sequence $(c_j)$ satisfying \ref{item_existence_of_c} can be defined by $c_j = \abs{\sigma_j}^\frac{\alpha}{2-\alpha}$.

Note that the conclusion in this example is not optimal in the case of $\alpha$-stable noise. For, in Example \ref{ex.one-sided-stable} we can choose $\sigma \in 
\ell^{\frac{2\alpha}{2-\alpha}}(\R)$ whereas here we have to choose
$\sigma\in \ell^{\frac{2\delta}{2-\delta}}(\R)$ for $\delta<\alpha$.
\end{Example}

The integration theory developed in \cite{Riedle_L2} and improved in Section \ref{sec_finite_second} relies on finite weak moments of the cylindrical L\'evy process. In the following, we extend this stochastic integral to the class of cylindrical L\'evy processes of the form \eqref{series_cylindrical_Levy} under Assumption \ref{item_existence_of_c} without requiring finite weak moments. For this purpose, by fixing a sequence $c\in \ell^\infty(\Rp)$ such that Assumption \ref{item_existence_of_c}  is satisfied and by using the notation \eqref{eq.def-p}--\eqref{eq.def-q} we define for each $k>0$:
\begin{equation*}
L^c_k(t)u := \sum_{j=1}^\infty \Big( p_j^{c,k}(t)+ m_j^{c,k}(t)\Big) \langle u, e_j \rangle, \qquad t\ge 0,\, u \in U.
\end{equation*}
Lemma  \ref{lem_cyl_decomp} yields that $L^c_k=P_k^c+M_k^c$ is a square-integrable cylindrical L\'evy process. 
Let $Q_k$ denote the covariance operator of $M_k^c$ and $\mathcal{H}_k$ be the corresponding reproducing kernel Hilbert space, where we suppress the dependence  on the sequence $c$ in the notation for $Q_k$ and $\mathcal{H}_k$.
At the same time, we extend the class of integrands by the usual localisation arguments. 
For this purpose, we define the class $\Lambda_{\rm{loc}}$ as
\begin{multline*}
\Lambda_{\rm{loc}} :=  \bigg\{ \Psi:[0,T]\times \Omega \to L_{\rm HS}(\mathcal{H}_k,H) : \Psi \text{ is predictable and } \\ \int_0^T \norm{\Psi(t)Q_k^{1/2}}_{L_{\rm HS}(U,H)}^2 \dt<\infty \quad P\text{-a.s. for all } k \in \N\bigg\}.
\end{multline*}

\begin{Theorem}\label{th.integration-predictable} 
Assume that $L$ is a cylindrical L\'evy process of the form \eqref{series_cylindrical_Levy} satisfying {\upshape\ref{item_existence_of_c}} 
for a sequence $c\in \ell^\infty(\Rp)$ and let $\Psi$ be in $\Lambda_{\rm{loc}}$. Then there exists an increasing sequence of stopping times $(\varrho(k))$  with $\varrho(k)\to \infty$ $P$-a.s.\ as $k\to\infty$ such that $\Psi(\cdot)\1_{[0,\varrho(k)]}(\cdot)\in \Lambda$ for each $k\in\N$ and 
\begin{align*}
\left(\int_0^t \Psi(s)\1_{\{s\leq\varrho(k)\}} \dL^c_k(s):\, t\in [0,T]\right)_{k\in\N}
\end{align*}
is a Cauchy sequence in the topology of uniform convergence in probability and its limit is independent of the sequence $c$ satisfying Assumption \ref{item_existence_of_c}.
\end{Theorem}

Theorem \ref{th.integration-predictable} enables us to define for each $\Psi\in \Lambda_{\rm loc}$ the stochastic integrals
\begin{align*}
\int_0^\cdot \Psi(s)\, \dL(s):=\lim_{k\to\infty} \int_0^\cdot \Psi(s)\1_{[0,\varrho(k)]}(s)\, \dL^c_k(s),
\end{align*}
where the limit is taken in the topology of uniform convergence in probability.
Note, that although in \cite{Jakubowski_Riedle} a stochastic integration theory is developed for a large class of integrands with respect to arbitrary cylindrical L\'evy processes, it does not cover the case of only predictable integrands.

\begin{proof}
Define the stopping times
$$\tilde{\tau}(k,n) := \inf \left\{ t\geq 0 : \int_0^t \bigHSnorm{\Psi(s)Q_k^{1/2}}^2 \ds > n \right\},$$
where we take the convention that $\inf \varnothing = + \infty$.
Since $\Psi \in \Lambda_{\rm{loc}}$, for every $k$ there is $n_k$ such that $P( \tilde{\tau}(k,n_k)<\infty)  \leq \frac{1}{2^k}$.
By the Borel--Cantelli Lemma
$$P\left( \limsup_{k\to \infty} \{ \tilde{\tau}(k,n_k)<\infty\} \right) = 0.$$
Consequently, the stopping times  $\varrho^c(k) := \tau^c(k) \wedge \tilde{\tau}(k,n_k) $ also converge to $+\infty$ a.s. by Lemma \ref{lem_cyl_decomp}. Note that if $T< \varrho^c(k)$, then $L^c_k=L^c_n$ on $[0,T]$ and
$$\int_0^t \Psi(s) \1_{\{s\leq\varrho^c(k)\}} \dL^c_k(s) = \int_0^t \Psi(s) \1_{\{s\leq\varrho^c(n)\}} \dL^c_n(s)$$
for all $t \in [0,T]$.
Consequently,  we obtain for each $k\le n$ and $\epsilon>0$ that 
\begin{align*}
&P\left( \sup_{t \in [0,T]} \normH{\int_0^t \Psi(s) \1_{\{s\leq\varrho^c(k)\}} \dL^c_k(s) - \int_0^t \Psi(s) \1_{\{s\leq\varrho^c(n)\}} \dL^c_n(s)} \ge \epsilon \right) \\
&\qquad \leq P \left( \int_0^t \Psi(s)\1_{\{s\leq\varrho^c(k)\}} \dL^c_k(s) \neq \int_0^t \Psi(s) \1_{\{s\leq\varrho^c(n)\}} \dL^c_n(s)\text{ for some } t \in [0,T] \right)\\
&\qquad\le P\big(T \geq \varrho^c(k)\big)
\to 0 \qquad \text{as } n,k\to \infty, 
\end{align*}
which establishes the claimed convergence. 
%Note, that under the assumptions of the theorem, the integral in \eqref{Cauchy_sequence_in_ucp} is defined as described in Section \ref{sec_finite_second}.

The limit of the Cauchy sequence does not depend on the choice of the sequence $c$ satisfying \ref{item_existence_of_c} because if $d$ is another sequence satisfying \ref{item_existence_of_c}, then  $L^c_k=L^d_n$ for all $t \in [0,T]$ on $\{T< \tau^c(k)\wedge \tau^d(n)\}$ and
$$\int_0^t \Psi(s) \1_{\{s\leq\tau^c(k)\}} \dL^c_k(s) = \int_0^t \Psi(s) \1_{\{s\leq\tau^d(n)\}} \dL^d_n(s),$$
which completes the proof.
\end{proof}

\section{Existence of a solution for the orthogonal noise}
\label{sec_existence}

Existence of a cylindrical L\'evy process of the form \eqref{series_cylindrical_Levy} strongly depends on the interplay between
the drift part $b_j$ and the L\'evy measure $\rho_j$ of the real valued L\'evy process with characteristics $(b_j,s_j,\rho_j)$, see condition \eqref{lemma_4_2_first_condition}. For this reason, we consider the general case of a cylindrical L\'evy process with a possible non-zero drift part. Naturally, we will tackle this part by moving it to the drift part of the equation under consideration. For this purpose, recall the decomposition 
$L(t)=P^c_k(t) +M^c_k(t)+R^c_k(t)$ of the cylindrical L\'evy process $L$ for each $k>0$ derived in Lemma \ref{lem_cyl_decomp} under assumption \ref{item_existence_of_c} satisfied for a sequence $c\in \ell^\infty(\Rp)$. 
Let $Q_k$ denote the covariance operator of $M_k^c$ where we suppress the dependency on the sequence $c$.
Furthermore, instead of the standard coercivity and monotonicity requirements, we introduce assumptions for each truncation level $k \in \N$. 
Assumptions of this form were introduced in Peszat and Zabczyk \cite[Sec.\ 9.7]{Peszat_Zabczyk} in the semigroup approach. Assume that there are constants $\alpha_k, \lambda_k, \beta_k>0$ such that
\begin{enumerate}[label=(A\arabic*$^\prime$)]
\item \label{item_monotonicty_not_L2}(coercivity) For every $k \in \N$ 
and $v\in V$ we have
\begin{align*}
2\pairing{F(v)+P_k^c(1) G^*(v)}v + \bigHSnorm{G(v)Q_{k}^{1/2}}^2 + \alpha_k \normV{v}^2 \leq \lambda_k \normH{v}^2+\beta_k;
\end{align*}
\item \label{item_coercivity_not_L2} (monotonicity) For every $k \in \N$ and $v_1,v_2 \in V$ we have
\begin{multline*}
2\pairing{F(v_1)-F(v_2)+P_k^c(1)\big(G^*(v_1)- G^*(v_2)\big)}{v_1-v_2} \\
+ \bigHSnorm{(G(v_1)-G(v_2))Q_{k}^{1/2}}^2 \leq \lambda_k \normH{v_1-v_2}. 
\end{multline*}
\item \label{item_linear_growth_not_L2} (linear growth) $\normVstar{F(v)+P_k^c(1)G^*(v)} \leq c_k(1+\normV{v} )$ for all $v \in V$;
\item \label{item_hemicontinuity_not_L2} (hemicontinuity) the mapping 
$$\R \ni s \mapsto {\pairing{F(v_1+sv_2) + P_k^c(1)G^*(v_1+sv_2)}{v_3}}$$
is continuous for all $v_1,v_2,v_3 \in V$.
\end{enumerate}

\begin{Theorem}
\label{th_existence_of_solution_not_L2}
Assume that $L$ is a cylindrical L\'evy process of the form \eqref{series_cylindrical_Levy} satisfying {\upshape \ref{item_existence_of_c}}  on page \pageref{item_existence_of_c} for a sequence $c\in \ell^\infty(\Rp)$.
If the coefficients $F$ and $G$ satisfy {\upshape (A1$^\prime$)}--{\upshape(A4$^\prime$)}, then equation \eqref{SPDE} with an $\mathcal{F}_0$-measurable initial condition $X(0) = X_0$ has a pathwise unique variational solution $(X,\bar{X})$. 
\end{Theorem}

\begin{proof}
We reduce the case of the general initial condition to the square integrable one as in \cite[Th.\ 6.2.3]{Applebaum}. For $k\in \N$ let $\Omega_k = \{ \norm{X_0} \leq k \}$ and $X_0^k = X_0 \1_{\Omega_k}$.
Using the decomposition $L(t)= P^c_k(t) +M^c_k(t)+R^c_k(t)$, Lemma \ref{lem_cyl_decomp} guarantees that $M^c_k$ is a weakly square-integrable cylindrical L\'evy martingale with diagonal covariance, and thus according to Theorem \ref{th_main_existence_result} there exists a unique variational solution $(X_k^c,\bar{X}_k^c)$ of 
%$$\dX(t) = F(X(t)) \dt + G(X(t)) \dL^k(t)$$
%so that it is driven by $\hat{L}^k$
$$\dX(t) = \big(F(X(t)) + P_k^c (1) G^*(X(t))) \dt + G(X(t)\big) \, \mathrm{d}M_k^c(t),$$
with the initial condition $X(0) = X_0^k$.

Step 1. 
We first show that for each $k\le n$ we have
$X_k^c(t)=X_n^c(t)$ $P$-a.s.\ on $\{t< \tau^c(k)\} \cap \Omega_k$.

Define a cylindrical L\'evy process $Y_{k,n}^c$ by 
$$Y_{k,n}^c(t)u := R_k^c(t)u - R_n^c(t)u = \sum_{j=1}^\infty \left( \int_{k/c_j < \abs{x} \leq n/c_j} x\,  N_j(t,\ud x)\right) \langle u,e_j \rangle$$
for all $t\geq 0$ and $u\in U$.
The cylindrical martingale $M_n^c$ can be rewritten as
\begin{align*}
M_n^c(t)u
&= \sum_{j=1}^\infty \left(m_j^{c,k}(t) + \int_{k/c_j < \abs{x} \leq n/c_j} x \, \tilde{N}_j(t,\mathrm{d}x) \right) \langle u,e_j \rangle  \\
&= \sum_{j=1}^\infty \left(m_j^{c,k}(t) + \int_{k/c_j < \abs{x} \leq n/c_j} x \, N_j(t,\ud x) - \int_{k/c_j < \abs{x} \leq n/c_j} x \, \rho_j(t,\mathrm{d}x) \right) \langle u,e_j \rangle \\
&=M_k^c(t)u+Y_{k,n}^c(t)u-(P_n^c(1)u-P_k^c(1)u)t.
\end{align*}
Applying this we get  
\begin{equation}
\label{difference_before_Ito_formula}
\begin{aligned}
X_k^c(t) - X_n^c(t)
&= -X_0\1_{\Omega_n \setminus \Omega_k} + \int_0^t \left( F\left(\bar{X}_k^c(s)\right)-F\left(\bar{X}_n^c(s)\right) \right) \ud s \\
&\quad + \int_0^t P_k^c(1) \left( G^*\left(\bar{X}_k^c(s)\right)-G^*\left(\bar{X}_n^c(s)\right) \right) \ud s \\ 
&\quad + \int_0^t \left( G\left(\bar{X}_k^c(s)\right) - G\left(\bar{X}_n^c(s)\right) \right) \ud M_k^c(s) \\
&\quad - \int_0^t G\left(\bar{X}_n^c(s)\right) \ud Y_{k,n}^c(s).
\end{aligned}
\end{equation}
By introducing the new notation
\begin{align*}
A(t) &:= X_k^c(t) - X_n^c(t) + \int_0^t G\left(\bar{X}_n^c(s)\right) \ud Y_{k,n}^c(s), \\
I(t) &:= \int_0^t  \left( G\left(\bar{X}_k^c(s)\right) - G\left(\bar{X}_n^c(s) \right) \right) \ud M_k^c(s),
\end{align*}
equation \eqref{difference_before_Ito_formula} can be rewritten as
\begin{equation}
\label{formula_for_A}
\begin{aligned}
A(t) 
&= -X_0\1_{\Omega_n \setminus \Omega_k} + \int_0^t \left( F\left(\bar{X}_k^c(s)\right)-F\left(\bar{X}_n^c(s)\right) \right) \ds \\
&\quad + \int_0^t P_k^c(1) \left( G^*\left(\bar{X}_k^c(s)\right)-G^*\left(\bar{X}_n^c(s)\right) \right) \ud s + I(t).
\end{aligned}
\end{equation}
Note that on $\{ t < \tau^c(k)\}$ we have $A(t)=X_k^c(t) - X_n^c(t)$.
We apply Theorem 1 in \cite{Gyongy_Krylov_II}  with $v(t) = \bar{X}_k^c(t) - \bar{X}_n^c(t)$ and $h(t) = I(t)$ on the set $\{ t < \tau^c(k)\}$.
% on which we have
%\begin{align*}
%X_k^c(t) - X_n^c(t)
%&= -X_0\1_{\Omega_n \setminus \Omega_k} + \int_0^t F\left(\bar{X}_k^c(s)\right)-F\left(\bar{X}_n^c(s)\right) \ds \\
%&\quad + \int_0^t P_k^c(1) \left( G^*\left(\bar{X}_k^c(s)\right)-G^*\left(\bar{X}_n^c(s)\right) \right) \ud s + I(t).
%\end{align*}
It follows that there exists an $H$-valued, c\`adl\`ag process $Z$, which is equal to $A=X_k^c - X_n^c$ $\Leb \otimes P$-almost everywhere on $\Upsilon := \{(t,\omega) \in [0,T] \times \Omega : t<\tau^c(k)(\omega) \}$ and such that the It\^{o} formula for the square of the norm holds:
\begin{equation}
\label{Ito_up_to_stopping_time}
\begin{aligned}
&\normH{Z(t)}^2 \\
&= \normH{X_0}^2 \1_{\Omega_n \setminus \Omega_k} + 2 \int_0^t \pairing{F\left(\bar{X}_k^c(s) \right)-F\left(\bar{X}_n^c(s)\right)}{\bar{X}_k^c(s) - \bar{X}_n^c(s)} \ds \\
&\quad + 2\int_0^t \pairing{P_k^c(1) \left( G^*\left(\bar{X}_k^c(s)\right)-G^*\left(\bar{X}_n^c(s)\right) \right)}{\bar{X}_k^c(s) - \bar{X}_n^c(s)} \ds \\ 
&\quad + \int_0^t \left( X_k^c(s-)-X_n^c(s-) \right) \ud I(s) + [I,I](t).\\
\end{aligned}
\end{equation}
We show that $X_k^c-X_n^c$ and $Z$ are indistinguishable on $\Upsilon$. 
We have that
$$\E \left[ \int_0^{\tau^c(k)} \1_{\{Z(t) \neq (X_k^c-X_n^c)(t)\}} \dt \right] = 0.$$
This implies that there exists $\Omega_1 \subset \Omega$ with $P(\Omega_1)=1$ such that for all $\omega \in \Omega_1$
$$\int_0^{\tau^c(k)(\omega)} \1_{\{Z(t,\omega) \neq (X_k^c-X_n^c)(t,\omega)\}} \dt = 0.$$
We obtain that for every $\omega \in \Omega_1$ there is a subset $A_{\omega,t} \subset [0,\tau^c(k)(\omega))$ with $\Leb(A_{\omega,t}) = \tau^c(k)(\omega)$ and such that $Z(t,\omega) = (X_k^c-X_n^c)(t,\omega)$ for all $t \in A_{\omega,t}$.
Note that $X_k^c-X_n^c$ is a c\`adl\`ag process in $V^*$. 
Fix $\omega \in \Omega_1$ and $t \in [0,\tau^c(k)(\omega))$. Let $(t_n) \subset A_{\omega,t}$ be a sequence decreasing to $t$. 
Since $Z(t_n,\omega) = (X_k^c-X_n^c)(t_n,\omega)$ for all $n\in \N$ we obtain that $Z(t,\omega) = (X_k^c-X_n^c)(t,\omega)$ for all $\omega \in \Omega_1$ and $t \in [0,\tau^c(\omega))$ i.e.\ $Z$ and $X_k^c-X_n^c$ are indistinguishable on $\Upsilon$.

Thus, in what follows, we can assume that for $t< \tau^c(k)$ the process $X_k^c-X_n^c$ is $H$-valued and c\`adl\`ag.
Let
$$J(t) := \int_0^t \left( X_k^c(s-)-X_n^c(s-) \right) \ud I(s).$$
We show that
\begin{equation}
\label{stopped_Ito_formula_to_prove}
\begin{aligned}
&\normH{A(t\wedge \tau^c(k))}^2 \\
&= \normH{X_0}^2 \1_{\Omega_n \setminus \Omega_k} + 2 \int_0^{t\wedge \tau^c(k)} \pairing{F\left(\bar{X}_k^c(s) \right)-F\left(\bar{X}_n^c(s)\right)}{\bar{X}_k^c(s) - \bar{X}_n^c(s)} \ds \\
&\quad + 2 \int_0^{t\wedge \tau^c(k)} \pairing{P_k^c(1) \left( G^*\left(\bar{X}_k^c(s)\right)-G^*\left(\bar{X}_n^c(s)\right) \right)}{\bar{X}_k^c(s) - \bar{X}_n^c(s)} \ds \\ 
&\quad + J(t\wedge \tau^c(k)) + [I,I](t \wedge \tau^c(k)).
\end{aligned}
\end{equation}
It follows from \eqref{Ito_up_to_stopping_time} by taking the left limit at $t \wedge \tau^c(k)$ that
\begin{equation}
\label{left_limit}
\begin{aligned}
&\normH{X_k^c((t\wedge \tau^c(k))-)-X_n^c((t\wedge \tau^c(k))-)}^2 \\
&= \normH{X_0}^2 \1_{\Omega_n \setminus \Omega_k} + 2 \int_0^{t\wedge \tau^c(k)} \pairing{F\left(\bar{X}_k^c(s) \right)-F\left(\bar{X}_n^c(s)\right)}{\bar{X}_k^c(s) - \bar{X}_n^c(s)} \ds \\
&\quad +2  \int_0^{t\wedge \tau^c(k)} \pairing{P_k^c(1) \left( G^*\left(\bar{X}_k^c(s)\right)-G^*\left(\bar{X}_n^c(s)\right) \right)}{\bar{X}_k^c(s) - \bar{X}_n^c(s)} \ds \\ 
&\quad + J \big( (t\wedge \tau^c(k))- \big) + [I,I] \big( (t \wedge \tau^c(k))- \big).
\end{aligned}
\end{equation}
By definition of $A$, $A(s)=X_k^c(s)-X_n^c(s)$ for $s<\tau^c(k)$. Taking the limits as $s\nearrow t \wedge \tau^c(k)$ we get that $A \big( (t\wedge \tau^c(k))- \big) = X_k^c \big( (t\wedge \tau^c(k))- \big) - X_n^c \big( (t\wedge \tau^c(k))- \big)$. 
Since the only discontinuous processes in \eqref{formula_for_A} are $A$ and $I$, it follows that $\Delta A(t\wedge \tau^c(k)) = \Delta I(t\wedge \tau^c(k))$.
Thus
\begin{align*}
\normH{A(t\wedge \tau^c(k))}^2 
&= \normH{A((t\wedge \tau^c(k))-) + \Delta A(t\wedge \tau^c(k))}^2 \\
&= \normH{X_k^c \big( (t\wedge \tau^c(k))- \big) - X_n^c \big( (t\wedge \tau^c(k))- \big)  + \Delta I(t\wedge \tau^c(k))}^2 \\
&= \normH{X_k^c \big( (t\wedge \tau^c(k))- \big) + X_n^c \big( (t\wedge \tau^c(k))- \big) }^2 + \normH{\Delta I(t\wedge \tau^c(k)}^2 \\
&\quad + \scalarH{\Delta I(t\wedge \tau^c(k))}{X_k^c \big( (t\wedge \tau^c(k))- \big) - X_n^c \big( (t\wedge \tau^c(k))- \big)}.
\end{align*}
Applying \eqref{left_limit} we obtain
\begin{equation}
\label{penultimate_form_of_A}
\begin{aligned}
&\normH{A(t\wedge \tau^c(k))}^2 \\
&= \normH{X_0}^2 \1_{\Omega_n \setminus \Omega_k} + \int_0^{t\wedge \tau^c(k)} \pairing{F\left(\bar{X}_k^c(s) \right)-F\left(\bar{X}_n^c(s)\right)}{\bar{X}_k^c(s) - \bar{X}_n^c(s)} \ds \\
&\quad + \int_0^{t\wedge \tau^c(k)} \pairing{P_k^c(1) \left( G^*\left(\bar{X}_k^c(s)\right)-G^*\left(\bar{X}_n^c(s)\right) \right)}{\bar{X}_k^c(s) - \bar{X}_n^c(s)} \ds \\ 
&\quad + J \big( (t\wedge \tau^c(k))- \big) + [I,I] \big( (t \wedge \tau^c(k))- \big) \\
&\quad + \scalarH{\Delta I(t\wedge \tau^c(k))}{X_k^c \big( (t\wedge \tau^c(k))- \big) - X_n^c \big((t\wedge \tau^c(k))- \big)} + \normH{\Delta I(t\wedge \tau^c(k))}^2
\end{aligned}
\end{equation}
The jump of the stochastic integral $J$ at $t\wedge \tau^c(k)$ equals to 
$$\scalarH{X_k^c((t\wedge \tau^c(k))-)-X_n^c((t\wedge \tau^c(k))-)}{\Delta I(t\wedge \tau^c(k))},$$
see \cite[Prop.\ 24.3 and Sec.\ 26.4]{Metivier}. Thus in \eqref{penultimate_form_of_A} we can simplify
\begin{multline}
\label{stochastic_integral_jump}
J \big( (t\wedge \tau^c(k))- \big) + \scalarH{\Delta I(t\wedge \tau^c(k))}{X_k^c \big( (t\wedge \tau^c(k))- \big)-X_n^c \big( (t\wedge \tau^c(k))- \big)} \\
\begin{aligned}
&= J \big(( t\wedge \tau^c(k))- \big) + \Delta J(t\wedge \tau^c(k)) \\
&= J(t\wedge \tau^c(k)).
\end{aligned}
\end{multline}
Similarly, the jump  of the quadratic variation of $I$ at $t \wedge \tau^c(k)$ 
equals $\normH{\Delta I(t \wedge \tau^c(k))}^2$, see \cite[Th.\ 20.5(4)]{Metivier}, and we can simplify in \eqref{penultimate_form_of_A}:
\begin{align}
\label{quadratic_variation_jump}
[I,I] \big( (t \wedge \tau^c(k))- \big) + \normH{\Delta I(t \wedge \tau^c(k \big))}^2 \notag
&= [I,I] \big((t \wedge \tau^c(k))-) +  \Delta [I,I](t \wedge \tau^c(k) \big) \\
&= [I,I](t \wedge \tau^c(k)).
\end{align}
Applying \eqref{stochastic_integral_jump} and \eqref{quadratic_variation_jump} in \eqref{penultimate_form_of_A} finishes the proof of \eqref{stopped_Ito_formula_to_prove}.

We multiply both sides of \eqref{stopped_Ito_formula_to_prove} by $\1_{\Omega_k}$ and take expectation. For the term involving the quadratic variation, we use the fact that $\E \left[ [I,I](t \wedge \tau^c(k) \right] = \E \left[ \langle I,I \rangle(t \wedge \tau^c(k) \right]$ and Theorem \ref{th_properties_integral}\ref{it_angle_bracket}.
Recall for the following that the martingale property is invariant under multiplication by $\1_{\Omega_k}$, since $\Omega_k$ is $\mathcal{F}_0$-measurable. 
Thus $\E \left[ J(t \wedge \tau^c(k))\1_{\Omega_k}  \right] = 0$. 
We obtain
\begin{equation*}
\begin{aligned}
&\E \left[ \normH{A(t\wedge \tau^c(k))}^2 \1_{\Omega_k} \right] \\
&= \E \left[ 2\int_0^{t\wedge \tau^c(k)} \pairing{F\left(\bar{X}_k^c(s) \right)-F\left(\bar{X}_n^c(s)\right)}{\bar{X}_k^c(s) - \bar{X}_n^c(s)} \ds \1_{\Omega_k} \right] \\
&\quad + \E \left[ 2\int_0^{t\wedge \tau^c(k)} \pairing{P_k^c(1) \left( G^*\left(\bar{X}_k^c(s)\right)-G^*\left(\bar{X}_n^c(s)\right) \right)}{\bar{X}_k^c(s) - \bar{X}_n^c(s)} \ds \1_{\Omega_k} \right] \\ 
&\quad + \E \left[ \int_0^{t\wedge \tau^c(k)} \bignorm{\left( G\left(\bar{X}_k^c(s)\right) - G\left(\bar{X}_n^c(s)\right) \right) Q_k^{1/2}}_{L_{\rm HS}(U,H)}^2 \ds \1_{\Omega_k} \right].
\end{aligned}
\end{equation*}
The monotonicity assumption \ref{item_monotonicty_not_L2} implies that
\begin{align*}
\E \left[ \normH{A(t\wedge \tau^c(k))}^2 \1_{\Omega_k} \right] 
&\leq \lambda_k \E \left[ \int_0^{t\wedge \tau^c(k)} \normH{X_k^c(s)-X_n^c(s)}^2 \ds \1_{\Omega_k} \right] \\
&\leq \lambda_k \E \left[ \int_0^t  \normH{A(s\wedge \tau^c(k))}^2 \ds \1_{\Omega_k} \right].
\end{align*}
It follows by Gronwall's inequality that 
$$\E \left[ \normH{X_k^c(t\wedge \tau^c(k))-X_n^c(t\wedge \tau^c(k)) + \int_0^{t\wedge \tau^c(k)} G\left(\bar{X}_n^c(s)\right) \ud Y_{k,n}^c(s)}^2 \1_{\Omega_k} \right] = 0.$$
Thus  
$$\left( X_k^c(t\wedge \tau^c(k))-X_n^c(t\wedge \tau^c(k)) + \int_0^{t\wedge \tau^c(k)} G\left(\bar{X}_n^c(s)\right) \ud Y_{k,n}^c(s) \right) \1_{\Omega_k} = 0 \quad \text{ a.s.}$$
In particular we obtain that 
$$X_k^c(t)-X_n^c(t) = 0 \text{ a.s. on } \{t<\tau^c(k)\} \cap \Omega_k.$$

Step 2. 
The first part enables us to define 
\begin{equation}
\label{def_of_X_by_X_k}
X := X^c_k \quad \text{ and } \quad \bar{X} := \bar{X}^c_k \qquad \text{on } \{t<\tau^c(k)\}\cap \Omega_k.
\end{equation}
This definition does not depend on the choice of the sequence $c$: for, if $d$ is another sequence in $\ell^\infty(\R)$ satisfying \ref{item_existence_of_c}, then one can show similarly as in Step 1, that $X^c_k=X^d_n$ on $\{ t < \tau^c(k) \wedge \tau^d(n)\}$. 
Moreover, since $X_k^c$ and $\bar{X}_k^c$ are equal almost everywhere on $\{ (\omega,t) \in \Omega_k \times [0,T] : t<\tau^c(k)(\omega)\}$, it follows by taking $k\to \infty$ that $X$ and $\bar{X}$ are equal almost everywhere on $\Omega \times [0,T]$. 

Step 3. We show that $(X,\bar{X})$ defined in \eqref{def_of_X_by_X_k} satisfies \eqref{integral_eq_in_def_of_variational_sol}. Note that
\begin{align}\label{eq.X-solution-on-tau}
&X(t) \1_{\{t<\tau^c(k)\}\cap \Omega_k} 
= X^c_k(t) \1_{\{t<\tau^c(k)\}\cap \Omega_k} \notag \\
& = X_0\1_{\Omega_k} + \1_{\{t<\tau^c(k)\}\cap \Omega_k}  \int_0^t F(\bar{X}^c_k(s)) \ds +\1_{\{t<\tau^c(k)\}\cap \Omega_k}  \int_0^t G(\bar{X}^c_k(s)) \dL^c_k(s) .
\end{align}
From the very definition \eqref{def_of_X_by_X_k} it follows 
\begin{align}\label{eq.F-and-Xck}
\lim_{k\to\infty} \1_{\{t<\tau^c(k)\}\cap \Omega_k} \int_0^t F(\bar{X}^c_k(s)) \ds 
&= \lim_{k\to\infty} \1_{\{t<\tau^c(k)\}\cap \Omega_k}  \int_0^t F(\bar{X}(s)) \ds\notag \\
&= \int_0^t F(\bar{X}(s)) \ds. 
\end{align}
The last term in \eqref{eq.X-solution-on-tau} can be rewritten as
\begin{align*}
\1_{\{t<\tau^c(k)\}\cap \Omega_k}  \int_0^t G(\bar{X}^c_k(s)) \dL^c_k(s)
=\1_{\{t<\tau^c(k)\}\cap \Omega_k}  \int_0^{t\wedge \tau^c(k)} G(\bar{X}^c_k(s)) \dL^c_k(s).
\end{align*}
From Theorem \ref{th_properties_integral}\ref{it_stable_under_stopping} and the definition of the stochastic integral with respect to $L$ after Theorem \ref{th.integration-predictable}, it follows that
\begin{align}\label{eq.G-and-Xck}
\lim_{k\to\infty}
 \int_0^{t\wedge \tau^c(k)} G(\bar{X}^c_k(s)) \dL^c_k(s)
&=\lim_{k\to\infty} \int_0^t G(\bar{X}^c_k(s))\1_{\{s \leq \tau^c(k)\}} \dL^c_k(s)\notag \\
&=\lim_{k\to\infty} \int_0^t G(\bar{X}(s))\1_{\{s\leq \tau^c(k)\}} \dL^c_k(s)\notag \\
&=\int_0^t G(\bar{X}(s)) \dL(s).
\end{align}
By taking the limit $k\to\infty$ in \eqref{eq.X-solution-on-tau}, 
equalities \eqref{eq.F-and-Xck} and \eqref{eq.G-and-Xck} show 
\begin{align*}
X(t)= X_0+\int_0^t F(\bar{X}(s))\ds \int_0^t G(\bar{X}(s)) \dL(s),
\end{align*}
which finishes the proof of the theorem.
\end{proof}

\subsection*{Acknowledgement}
We would like to thank a referee of an earlier version of this paper who pointed out some inconsistencies in Section \ref{sec_finite_second}, and which resulted in this improved version of the article.

%\bibliography{References_variational}{}
\bibliographystyle{plain}

\end{document}